\documentclass[11pt]{article}
\usepackage{amssymb, amsmath, url, graphicx, setspace, geometry}
\usepackage{calrsfs}
\usepackage{wasysym}
\usepackage{xcolor}

\def\3{\subset }
\def\4{\subseteq }
\def\<{\left<}
\def\>{\right>}

\def\bit{\begin{itemize}}
\def\eit{\end{itemize}}
\def\3{\subset }
\def\4{\subseteq }

\def\0{\leqno}

\def\barr{\begin{array}}
\def\earr{\end{array}}

\def\Z{{\rlap{$\kern2pt{\rm Z}$}{\rm Z}\,}}

\title{\bf Addendum to ``A generalization of a result on the sum of element orders of a finite group"}
\author{Mihai-Silviu Lazorec and Marius T\u arn\u auceanu}

\begin{document}
\maketitle

\begin{abstract}
Let $G$ be a group of order $n$ and $H$ be a subgroup of order $m$ of $G$. Denote by $\psi_H(G)$ the sum of element orders relative to $H$ of $G$. It is known that if $G$ is nilpotent, then $\psi_H(G)\leq\psi_{H_m}(G)$, where $H_m$ is the unique subgroup of order $m$ of $C_n$. In this note, we show that this inequality does not hold for infinitely many finite solvable groups.
\end{abstract}



\section{Introduction} 

Let $G$ be a finite group and $H$ be a subgroup of $G$. We denote by $o_H(x)$ the order of an element $x\in G$ relative to $H$, i.e. the smallest positive integer $m$ such that $x^m\in H$. Also, $C_n$ denotes the cyclic group with $n$ elements, where $n$ is a positive integer. 

The main objective of \cite{9} was to study the quantity 
$$\psi_H(G)=\sum\limits_{x\in G}o_H(x),$$
which is a generalization of the sum of element orders of $G$, that is
$$\psi(G)=\sum\limits_{x\in G}o(x),$$
where $o(x)$ denotes the usual order of an element $x\in G$. The second quantity appeared in \cite{1} and was further studied in several papers such as \cite{2,3,4,8}. The main result of \cite{1} states that the maximum value of the sum of element orders among all finite groups of order $n$ is $\psi(C_n)$. In other words, we have:\\

\textbf{Theorem 1.1.} \textit{Let $G$ be a group of order $n$. Then $\psi(G)\leq \psi(C_n)$ and the equality holds if and only if $G\cong C_n$.}\\ 

Up to the class of finite nilpotent groups, in \cite{9}, it was proved that a similar property occurs if one works with element orders relative to a subgroup. More exactly, the following result holds:\\

\textbf{Theorem 1.2.} \textit{Let $G$ be a nilpotent group of order $n$ and $H$ be a subgroup of order $m$ of $G$. Then 
\begin{align}\label{r1}
\psi_H(G)\leq \psi_{H_m}(C_n),
\end{align}
where $H_m$ is the unique subgroup of order $m$ of $C_n$.}\\

Moreover, the author conjectured that Theorem 1.2 is valid for any finite group. However, in this paper we show that this statement is false by outlining a way to construct a finite group $G$ such that inequality (\ref{r1}) is not satisfied for at least one subgroup $H$ of $G$. For the ease of writing, we introduce the following ratio
$$\psi'_H(G)=\frac{\psi_H(G)}{\psi_{C_{|H|}}(C_{|G|})}.$$
Using this notation, our main result may be stated as follows:\\

\textbf{Theorem 1.3.} \textit{There are infinitely many finite groups $G$ having a subgroup $H$ such that $\psi'_{H}(G)>1$.}\\

A proof of Theorem 1.3 along with other several properties related to this generalization of the sum of element orders of a finite group are highlighted in the following section. 

\section{Proof of the main result and new properties of $\psi_H(G)$}

As a starting point, we outline a preliminary result which puts together several known facts concerning the sum of element orders of a finite group $G$ and its generalization. For an element $x\in G\setminus H$, a novelty would be an upper bound for $o_H(x)$. As a consequence, we also determine a new upper bound for $\psi_H(G)$.\\

\textbf{Lemma 2.1.} \textit{Let $G$ be a group of order $n=p_1^{\alpha_1}p_2^{\alpha_2}\ldots p_k^{\alpha_k}$, where $p_1<p_2<\ldots<p_k$ are primes and $\alpha_i$ is a positive integer, for all $i\in\lbrace 1, 2,\ldots, k\rbrace$, $H$ be a subgroup of $G$ of index $q$, $K$ be a normal subgroup of $G$ and $(G_i)_{i=\overline{1,s}}$ be a family of finite groups of coprime orders. Then:
\begin{itemize}
\item[i)] $\psi(C_n)=\prod\limits_{i=1}^k\frac{p_i^{2\alpha_i+1}+1}{p_i+1};$ (see Lemma 2.9 (1) of \cite{4})
\item[ii)] $\psi(C_n)\geq\frac{p_1}{p_k+1}n^2;$ (see Lemma 2.9 (2) of \cite{4})
\item[iii)] $\psi_K(G)=|K|\psi(\frac{G}{K})$; (see p. 628 of \cite{9})
\item[iv)] $\psi_{H_1\times H_2\times\ldots\times H_s}(G_1\times G_2\times\ldots\times G_s)=\prod\limits_{i=1}^s\psi_{H_i}(G_i)$, where $H_i$ is a subgroup of $G_i$, for all $i\in\lbrace 1,2,\ldots, s\rbrace$; (see Lemma 2.1 a) of \cite{9})
\item[v)] $o_H(x)\leq q, \ \forall \ x\in G\setminus H;$
\item[vi)] $\psi_H(G)\leq |H|(q^2-q+1).$
\end{itemize}}
\textbf{Proof.} For item \textit{v)}, we have
$$o_H(x)=\frac{|\langle x\rangle|}{|\langle x\rangle\cap H|}=\frac{|H\langle x\rangle|}{|H|}\leq \frac{|G|}{|H|}=q, \ \forall \ x\in G\setminus H.$$

The above inequality is further used to prove the result given by item \textit{vi)}. More exactly, one can write
$$\psi_H(G)=|H|+\sum\limits_{x\in G\setminus H}o_H(x)\leq |H|+q(|G|-|H|)=|H|(q^2-q+1).$$ 
\hfill\rule{1,5mm}{1,5mm}\\

Before proving our main result, we recall that a subgroup $H$ of a finite group $G$ is isolated in $G$ if, for every element $x\in G$, we have either $x\in H$ or $\langle x\rangle\cap H=1$. For more details about these subgroups, we refer the reader to \cite{5}. Also, in \cite{10}, the isolatedness property was characterized as follows: $H$ is isolated in $G$ if and only if $\psi_H(G)=|H|+\psi(G)-\psi(H).$ \\

\textbf{Proof of Theorem 1.3.} Let $p$ be a prime and $r\geq 3$ be an  integer. By following Example 7 of \cite{7}, starting with the field $\mathbb{F}_q$ with $q=p^r$ elements, we can construct a Frobenius group $\overline{G}=\overline{K}\rtimes \overline{H}$, whose kernel $\overline{K}$ is isomorphic to the additive group of $\mathbb{F}_q$, while $\overline{H}$ is isomorphic to the multiplicative group of $\mathbb{F}_q$. Hence, $\overline{K}\cong C_p^r$ and $\overline{H}\cong C_{p^r-1}$. Moreover, since $\overline{H}$ is isolated in $\overline{G}$ and its conjugates trivially intersect, we obtain
$$o_{\overline{H}}(\overline{G})=|\overline{H}|+\psi(\overline{G})-\psi(\overline{H})=(p^r-1)(\psi(\overline{H})+p)=(p^r-1)(\psi(C_{p^r-1})+p).$$
Further, suppose that $p=2$ and $r$ is chosen such that $2^r-1$ is a prime (i.e. $2^r-1$ is a Mersenne prime). Using Lemma 2.1 \textit{i)}, \textit{iii)}, we get
$$\psi'_{\overline{H}}(\overline{G})=\frac{o_{\overline{H}}(\overline{G})}{o_{C_{|\overline{H}|}}(C_{|\overline{G}|})}=\frac{(2^r-1)(\psi(C_{2^r-1})+2)}{(2^r-1)\psi(C_{2^r})}
=\frac{3\cdot 2^{2r}-9\cdot 2^r+15}{2^{2r+1}+1}>1.$$
Consequently, if $r\geq 3$ is chosen such that $2^r-1$ is a Mersenne prime, the Frobenius group $\overline{G}$ is an example of a finite group for which inequality (\ref{r1}) does not hold. Finally, to obtain infinitely many such examples, take $G=\overline{G}\times C_q$, where $q$ is any odd prime such that $q\nmid 2^r-1$. Then, $H=\overline{H}\times C_q$ is a subgroup of $G$. By Lemma 2.1 \textit{iv)}, we deduce that
$$\psi'_H(G)=\psi'_{\overline{H}}(\overline{G})\psi'_{C_q}(C_q)=\psi'_{\overline{H}}(\overline{G})>1,$$
and this concludes the proof.   
\hfill\rule{1,5mm}{1,5mm}\\

Using the notations from our previous proof, note that $\psi'_{\overline{H}}(\overline{G})=\frac{3\cdot 2^{2r}-9\cdot 2^r+15}{2^{2r+1}+1}$ is an increasing function with respect to $r$ and $\displaystyle \lim_{r \to\infty}\psi'_{\overline{H}}(\overline{G})=\frac{3}{2}$. Hence, the smallest of our examples of finite groups for which inequality (\ref{r1}) does not hold is obtained for $p=2$ and $r=3$. In this case, we have $\overline{G}\cong C_2^3\rtimes C_7$ (SmallGroup(56,11)) and $\psi'_{\overline{H}}(\overline{G})=\frac{45}{43}.$ Also, note that all the groups that were outlined in our previous proof are solvable so, in general, Theorem 1.2 does not hold for this class of finite groups. 

An unsolved problem in group theory is the following one: \textit{Given a group $G$ of order $n$, does there necessarily exist a bijection $f$ from $G$ onto $C_n$ such that for each element $x\in G$, we have $o(x)|o(f(x))?$} (see Problem 18.1 posed by I. M. Isaacs in \cite{6}). The answer is true if $G$ is a solvable group. We can formulate a similar question concerning element orders relative to subgroups: \textit{Let $G$ be a group of order $n$, $H$ be a subgroup of order $m$ of $G$ and $H_m$ be the unique subgroup of order $m$ of $C_n$. Does there necessarily exist a bijection $f$ from $G$ onto $C_n$ such that for each element $x\in G$, we have $o_H(x)|o_{H_m}(f(x))?$} As a consequence of the proof of Theorem 1.3, we can state that the answer is negative even for solvable groups.



In \cite{9}, it was proved that inequality (\ref{r1}) holds whenever $H$ is a normal subgroup of a finite group $G$. In this paper, we show that this property also occurs for any subgroup $H$ of $G$ of prime index.\\

\textbf{Proposition 2.2.} \textit{Let $G$ be a group of order $n=mq$ and $H$ be a subgroup of $G$ of prime index $q$. Then $\psi_{H}(G)\leq \psi_{H_m}(C_n)$, where $H_m$ is the unique subgroup of order $m$ of $C_n$.}

\textbf{Proof.} By Lemma 2.1 \textit{i)}, \textit{iii)}, we have
$$\psi_{H_m}(C_n)=|H_m|\psi(C_q)=m(q^2-q+1).$$
Then, the result easily follows by Lemma 2.1 \textit{vi)}.
\hfill\rule{1,5mm}{1,5mm}\\

Finally, we obtain some upper bounds for $\psi'_H(G)$ in terms of the prime divisors of the index of $H$ in $G$.\\

\textbf{Proposition 2.3.} \textit{Let $G$ be a finite group and $H$ be a subgroup of $G$ of index $q=p_1^{\alpha_1}p_2^{\alpha_2}\ldots p_k^{\alpha_k}$, where $p_1<p_2<\ldots<p_k$ are primes and $\alpha_i$ is a positive integer, for all $i\in\lbrace 1, 2,\ldots, k\rbrace$. Then,
\begin{itemize}
\item[i)] $\psi'_H(G)<\prod\limits_{i=1}^k\frac{p_i+1}{p_i}\leq (\frac{3}{2})^k$;
\item[ii)] $\psi'_H(G)<\frac{p_k+1}{p_k}$; in particular, if $k=1$, then $\psi'_H(G)<\frac{p_1+1}{p_1}\leq\frac{3}{2}.$ 
\end{itemize}}      
\textbf{Proof.} By Lemma 2.1 \textit{iii)}, \textit{vi)}, we get
$$\psi'_H(G)=\frac{\psi_H(G)}{\psi_{C_{|H|}}(C_{|G|})}\leq \frac{q^2-q+1}{\psi(C_q)}.$$
Consequently, for item \textit{i)}, we apply Lemma 2.1 \textit{i)} to obtain
$$\psi'_H(G)\leq(q^2-q+1)\prod\limits_{i=1}^k\frac{p_i+1}{p_i^{2\alpha_i+1}+1}<\frac{q^2-q+1}{q^2}\prod\limits_{i=1}^k\frac{p_i+1}{p_i}<\prod\limits_{i=1}^k\frac{p_i+1}{p_i}\leq \bigg (\frac{3}{2}\bigg )^k,$$
while, for item \textit{ii)}, we use Lemma 2.1 \textit{ii)} to get
$$\psi'_H(G)\leq \frac{q^2-q+1}{q^2}\cdot \frac{p_k+1}{p_1}<\frac{p_k+1}{p_1}.$$
\hfill\rule{1,5mm}{1,5mm}\\

In what concerns the last result and its proof, note that it is not always true that $f(q)=\frac{q^2-q+1}{\psi(C_q)}<\frac{3}{2}$. For instance, if we take $q=2^a\cdot 3$, where $a\geq 1$ is an integer, we have
$$\displaystyle \lim_{a \to\infty}f(q)=\displaystyle \lim_{a \to\infty}\frac{9\cdot 2^{2a}-3\cdot 2^a+1}{\frac{7}{3}(2^{2a+1}+1)}=\frac{27}{14}>\frac{3}{2}.$$ 

\bigskip\noindent {\bf Acknowledgements.} This work was supported by a grant of the "Alexandru Ioan Cuza" University of Iasi, within the Research Grants program, Grant UAIC, code GI-UAIC-2021-01.

\vspace*{3ex}
\small

\begin{minipage}[t]{7cm}
Mihai-Silviu Lazorec \\
Faculty of  Mathematics \\
"Al.I. Cuza" University \\
Ia\c si, Romania \\
e-mail: {\tt silviu.lazorec@uaic.ro}
\end{minipage}
\hspace{3cm}
\begin{minipage}[t]{7cm}
Marius T\u arn\u auceanu \\
Faculty of  Mathematics \\
"Al.I. Cuza" University \\
Ia\c si, Romania \\
e-mail: {\tt tarnauc@uaic.ro}
\end{minipage}
\end{document}